\documentclass[twoside,10pt,reqno]{amsart}

\usepackage{upref,amsxtra,amssymb,amscd}
\usepackage{varioref}
\usepackage{verbatim}
\usepackage{epsfig}

\usepackage[french,english]{babel}
\usepackage{epic,eepic,eucal}
\usepackage[T1]{fontenc}
\usepackage{enumerate}

\def\jk{ \mathbb S}

\def\a{         \alpha}
\def\T{ {\bf T}   }
\def\R{ {\bf R}   }

\newcommand{\NN}{{\mathbb N}}
\newcommand{\RR}{{\mathbb R}}
\newcommand{\TT}{{\mathbb T}}
\newcommand{\ZZ}{{\mathbb Z}}
\newcommand{\CC}{{\mathbb C}}
\newcommand{\QQ}{{\mathbb Q}}

\def\carre{ \hfill $\Box$    }

\newtheorem*{prop*}{\sc Proposition}
\newtheorem*{lemm*}{\sc Lemma}
\newtheorem*{coro*}{\sc Corollary}

\newtheorem*{theo*}{\sc Theorem}
 
\theoremstyle{definition}

\newtheorem*{defi*}{\sc Definition}

\theoremstyle{remark}

\newtheorem{rema}{\sc Remark}

\numberwithin{equation}{section}

 \def\tn{|\kern-1pt|\kern-1pt|}

\begin{document}
 \title{Mixing diffeomorphisms and flows with purely singular spectra.}
\author{Bassam R. Fayad}

\maketitle

\begin{abstract} We give a geometric criterion that guaranteesa purely singular  spectral type for a dynamical system on a Riemannian manifold. The criterion, that is based on the existence of fairly rich but localized periodic approximations, is compatible with mixing. Indeed, we use it to construct examples of smooth mixing flows on the three torus with purely singular spectra. \end{abstract}

\maketitle

\section{Introduction} 

\subsection{}  \label{application} Mixing is one of the principal characterisctics of stochastic 
behavior of a dynamical system $(T,M,\mu)$. It is a spectral property and in the 
great majority of cases it is  a consequence of much stronger 
properties  of the system,  such as the K-property or fast correlation decay, which imply a Lebesgue spectrum for the associated unitary operator: $f \rightarrow f \circ T$ defined on $L^2(M,\mu,\CC)$.  In this paper, we prove the following 

\begin{theo*} There exist on $\TT^d$, $d \geq 3$, volume preserving flows and diffeomorphisms of class $C^\infty$ 
that are mixing and have purely singular spectra. 
\end{theo*} 

The only previously known examples 
where mixing of the system was accompanied by singular spectrum of the associated unitary operator were obtained in an abstract measure theoretical or probabilistic frame, such as Gaussian and 
related systems which by their nature do not come from 
differentiable dynamics, or rank one and mixing constructions which do not have yet $C^\infty$ realizations. 
  
To obtain our examples, we introduce a criterion for singular spectrum, based on the existence of fairly rich families of almost periodic sets, that is compatible with mixing, albeit at a slow rate. Then, we construct smooth mixing reparametrizations of some Liouvillean linear flows on $\TT^3$ satisfying the criterion, hence displaying a purely singular spectrum. 
As a by-product we observe,   due to Host's theorem \cite{host}, that the latter mixing reparametrizations are actually mixing of all orders.

\medskip

\subsection{Periodic approximations and singular spectra.} \label{PASS} $ \ $

\vspace{0.2cm}

A basic property implying  the singularity of the spectrum of $(T,M,\mu)$ is {\it rigidity}, i.e. the existence of a sequence of times $t_n $ such that for any measurable set $A \subset M$ it holds that $\mu \left(T^{t_n} A \triangle A \rightarrow 0 \right)$ where $A \triangle B$ stands for the symetric difference between $A$ and $B$. For smooth systems the latter property is often obtained as a consequence of a stronger one, namely the existence of cyclic (or more generally periodic) approximations in the sense of Katok and Stepin, see \cite{Kams}, i.e. the existence of a sequence of almost periodic towers such that any measurable set can  asymptotically be approximated by levels from the indivual towers. 

Rigidity of $(T,M,\mu)$ is clearly not compatible with mixing. To get a criterion that guarantees a singular spectrum without precluding mixing, we relax the concept of periodic approximations to that of having
strongly periodic towers with nice levels whose total measure might tend to zero but such that any measurable set can be covered by unions of levels from possibly different towers.

 \begin{defi*}[Slowly coalescent periodic approximations]  Let $T$ be an ergodic transformation of a Riemannian manifold $M$ preserving a volume $\mu$. We say that the dynamical system $(T,M, \mu)$ dispalys {\sl slowly coalescent  periodic approximations}, if there exists $\gamma >1$ and a sequence of integers $k_{n+1} \geq \gamma^n k_n$ such that for every $n \in \NN$ there exists a sequence 
$${\mathcal C}_n =  \bigcup_{i \in \NN} B_{n,i}$$
where the $B_{n,i}$, $i=0, \dots$, are {\sl balls} of $M$ satisfying
\begin{enumerate}
\item[(i)] $\displaystyle{\sup_{i \in \NN } r(B_{n,i}) \mathop{\longrightarrow} \limits_{n \rightarrow \infty} 0}$,
\item[(ii)] $\displaystyle{  \mu \left( T^{k_n} B_{n,i} \triangle B_{n,i} \right) \leq {\gamma}^{-n} \mu(B_{n,i}),}$ (where $\triangle$ denotes the symetric difference between sets),
\item[(iii)] $\mu \left( \mathop{\bigcap} \limits_{m \in \NN} \mathop{\bigcup} \limits_{n \geq m} {\mathcal C}_n \right) = 1$.
\end{enumerate}
\end{defi*} 

In Section \ref{S2} we will prove the following theorem

\begin{theo*}[Criterion for the singularity of the spectrum]  
A dynamical system $(T,M, \mu)$ dispalying {\sl slowly coalescent  periodic approximations} has a purely singular spectral type.
\end{theo*}

\begin{rema} In general, $\mu ({\mathcal C}_n)$ need not converge to zero. For a rotation of the circle, for example, it tends to the contrary to 1, in which case the terminology {\sl slowly coalescent} becomes a euphemism. For a mixing system $(T, M, \mu)$ however, $(ii)$ implies that $\mu ({\mathcal C}_n) \rightarrow 0$ and this is what we refer to by {\sl coalescent}. The terminology {\sl slowly coalescent}  is then used to refer to property $(iii)$ that is the key property in guaranteeing a purely singular spectrum. We will abreviate   slowly coalescent periodic approximations with SCPA.
\end{rema} 

\begin{rema} If the sets ${\mathcal C}_n$ satisfy adequate independence conditions, $(iii)$ will follow from the Borel Cantelli Lemma if $\sum \mu( {\mathcal C}_n ) = + \infty$.  \end{rema} 


\medskip

\subsection{Spectral type of reparametrized linear flows.} \label{repara} $ \ $ 

\vspace{0.2cm}

The problem  of understanding the ergodic 
and spectral properties of reparametrizations of  linear flows
on tori were raised by A. N. Kolmogorov in his I.C.M. address of 1954 
\cite{Ko}.  Since then and starting with the work of Kolmogorov 
himself,  this problem has been intesively studied and a surprisingly 
rich variety of behaviors were discovered  to be possible for the 
reparametrized flows. We say surprisingly because at the time when 
Kolmogorov raised the problem, some strong restrictions on the  spectral type of the reparametrized flow were expected to hold, at least in the 
case of real analytic reparametrizations, Cf. \cite{Ko} as well as the appendix by Fomin to the russian version of the book of Halmos on ergodic theory where absence of mixed spectrum was conjectured for smooth reparametrizations of linear flows.

We denote by $R^t_\alpha$ the linear flow on the torus $\TT^{n}$  given by
$$ \frac{d x}{d t} = \alpha, $$
where $x \in \TT^{n}$ and $\alpha$ is a vector of ${\RR}^n$. Given 
a continuous function $\phi : \TT^{n} \rightarrow {\RR}^{*}_{+}$ we
define the reparameterization flow $T^t_{\alpha, \phi}$ by
$$  \frac{d x}{d t} = \frac{ \alpha}{\phi (x,y)}. $$

  If the coordinates of $\alpha$ are rationally independent then the 
linear flow $R^t_{\alpha}$ is uniquely ergodic and so is $T^t_{\a,\phi}$ that preserves the measure with density $\phi$. 
Other properties of the linear flow may change under 
reparameterization.  While the linear flow has discrete (pure point) 
spectrum with the group of eigenvalues isomorphic to ${\ZZ}^{n}$, a continuous time change may yield a wide
variety of spectral properties. This follows from the theory of 
monotone (Kakutani) equivalence \cite{K} and the fact that every 
monotone measurable time change is cohomologous to a continuous one 
\cite{OS}. However, for 
sufficiently smooth reparameterizations the possibilities are more 
limited and they depend on the arithmetic properties of the vector 
$\alpha$.

If $\alpha$ is Diophantine and the function $\phi$ is $C^{\infty}$, 
then the reparameterized flow is smoothly isomorphic to a linear 
flow. This was first noticed by A.~N. Kolmogorov \cite{Ko}.  Herman found in 
\cite{He} sharp results of that kind for the finite regularity 
case. Kolmogorov also knew that for a Liouville vector $\alpha$ a 
smooth reparametrization could be weak mixing, or equivalently the associated unitary operator to the reparametrized flow could have 
a continuous spectrum.

 M.~D. {\v{S}}klover 
 proved in \cite{Shklover} the existence of real-analytic weak mixing 
reparametrizations of some Liouvillean linear flows on $\TT^2$; his 
result being optimal in that
he showed that for any real-analytic reparametrization $\phi$ other than a 
trigonometric polynomial there is $\alpha$ such that $T^t_{\alpha, 
\phi}$ is weakly mixing. In  \cite{weakmixing}, it was shown that  for any Liouvillean translation 
flow $R^t_\alpha$ on the torus
$\TT^n$, $n \geq 2$, the  generic $C^\infty$ reparametrization of 
$R^t_\alpha$ is weakly mixing.

Continuous and discrete spectra are not the only possibilities. In 
\cite{FKW}, B. Fayad, A. Katok and A. Windsor have proved that for 
every  $\alpha \in {\RR}^2$ with a Liouvillean slope there exists a 
strictly positive $C^\infty$ function
$\phi$ such that the flow on ${\TT}^2$ $T^t_{\alpha, \phi}$ has a 
mixed spectrum since it has a discrete part generated by only one 
eigenvalue. They also construct real-analytic examples for a more 
restricted class of Liouvillean $\alpha$. 

Recently, M. Guenais and F. Parreau  \cite{GP} achieved real-analytic 
reparametrizations of linear flows on ${\TT}^2$ that have an 
arbitrary number of eigenvalues. They even construct an example of a 
reparametrization of a linear flow on ${\TT}^2$ that is isomorphic to 
a linear flow on ${\TT}^2$ with "exotic" eigenvalues, i.e. not in the 
span of the eigenvalues of the original linear flow.

Finally, unlike continuous or descrete spectra, there exist real-analytic functions $\phi$ that are not  trigonometric polynomials, and for which a 
mixed spectrum is precluded for the flow $T^t_{\alpha, \phi}$  for any choice of $\alpha$.
  Indeed, it was proven in \cite{Dichotomy}  that
for a class of functions satisfying some regularity conditions on 
their Fourier coefficients the following dichotomy holds: 
$T^t_{\alpha, \phi}$ either has a continuous spectrum or is $L^2$ 
isomorphic to a constant time suspension.

\vspace{0.2cm} 

\noindent {\sl Reparametrizations and mixing.}  Katok \cite{K1} showed 
that for a function $\phi > 0$ of class
  $C^5$ any reparametrized flow $T^t_{\alpha, \phi}$ has
a simple spectrum, a singular maximal spectral type, and cannot be mixing. The singularity of the spectrum  was extended by A.~V. Ko\v{c}ergin to Lipschitz 
reparametrizations \cite{Kochergin1}.  The argument is based on a 
Denjoy--Koksma type estimate which fails in higher dimension 
\cite{Y}.  Based on the latter fact, it was  shown in \cite{mixing} that there exist $\alpha \in 
\RR^3$ and a real-analytic strictly positive function $\phi$ defined on $\TT^3$, such that the reparametrized 
flow $T^t_{\alpha, \phi}$  is mixing.
 
Recently Ko\v{c}ergin showed that for H\"older reparametrizations of 
some Diophantine linear flows on $\TT^2$ mixing is  possible \cite{Kodiophantine}.


The mixing examples obtained by reparametrizations of linear flows 
belong to a variety of fairly slow mixing systems, also including the 
mixing flows on surfaces constructed by Ko\v{c}ergin in the seventies 
\cite{Kochergin2}, for which the type and the multiplicity  of the spectrum 
remain undetermined. 

Modifying the reparametrizations of \cite{mixing}, it 
is possible to maintain mixing while the time one map of the reparametrized flow is forced to satisfy the SCPA criterion stated above, thus yielding

\begin{theo*} For $d \geq 3$, there exists $\alpha \in 
\RR^d$ and a strictly positive function $\phi$ over $\TT^d$ of class $C^\infty$ such that the reparametrized 
flow $T^t_{\alpha, \phi}$ is mixing and has a singular maximal spectral type with respect to the Lebesgue measure.  \end{theo*}

A dynamical system $(T,M,\mu)$ (or flow $(T^t,M,\mu)$) is said to be mixing of  order  $l \geq 2$ if, for any sequence $\displaystyle{ {\lbrace (u_n^{(1)},\cdots, u_n^{(l-1)}) \rbrace}_{n \in \NN}}$, where for $i=1, \cdots,l-1$ the ${\lbrace u_n^{(i)} \rbrace}_{n \in \NN}$ are sequences of integers (or real numbers) such that $\displaystyle{ \lim_{n \rightarrow  \infty} u_n^{(i)} = \infty}$, and for any $l$-upple $(A_1,\cdots,A_{l})$ of measurable subsets of $M$, we have 
$$ \lim_{n\rightarrow \infty} \mu \left(T^{-u_n^{(1)} - \ldots - u_n^{(l-1)}} A_{l} \cap \cdots \cap T^{-u_n^{(1)}} A_2 \cap A_1 \right) = \mu(A_{l-1}) \cdots \mu(A_1).$$

The general definition of mixing corresponds to mixing of order 2. A system is said to be mixing of all orders if it is mixing of order $l$ for any $l \geq 2$. Host's theorem \cite{host} asserts that a mixing system with singular spectrum is mixing of all orders, hence we get

\begin{coro*} 
  For $d \geq 3$, there exists $\alpha \in 
\RR^d$ and a strictly positive function $\phi$ over $\TT^d$ of class $C^\infty$ such that the reparametrized 
flow $T^t_{\alpha, \phi}$ is mixing of all orders.
\end{coro*}


\vspace{0.2cm} 

The paper consists of two sections. In Section \ref{S2} we prove Theorem-Criterion \ref{PASS}. In Section \ref{S3} we apply the criterion to obtain Theorem \ref{repara}.

 \medskip

\section{Slowly coalescent  periodic approximations} \label{S2}

In this section we prove Theorem \ref{PASS}. 
\subsection{} \label{spectral} We will use the following criterion that guarantees a singular spectrum for $(T,M,\mu)$:
\begin{prop*}  Let $(T, M, \mu)$ be a dynamcial system. If for any complex nonzero function $f \in L^2_0 (M, \mu)$, i.e. $\int_M f(x) d\mu (x)=0$, there exists a measurable set $E \subset M$ with $\mu (E) > 0$, and a strictly increasing sequence $l_n$, such that for every $x \in E$ we have 
\begin{eqnarray} \label{eq1} \limsup_{n \rightarrow \infty} {1 \over n} \left| \sum_{i=0}^{n-1} f(T^{l_i}x) \right| >0
\end{eqnarray}
then the maximal spectral type of the unitary operator associated to $(T,M, \mu)$ is singular.
\end{prop*} 

\noindent {\sl Proof.} Assume that $T$ has an absolutely continuous component in its spectrum. Then there exists $f \in L^2_0(M,\mu)$ such that the spectral measure corresponding to $f$ on the circle $\jk$ writes as $\sigma_f(dx) = g(x) dx$ where $g \in L^1(\jk, \RR_+, dx)$ is bounded. 
With the notation 
$$S_nf(x) = \sum_{i=0}^{n-1} f(T^{l_i}(x)) $$
we write spectrally
\begin{eqnarray*} 
{\left\| {S_nf  \over n} \right\|}_{L^2}^2 &=& {1 \over n^2} { \int_{\jk}}  {\left| \sum_{i=0}^{n-1} z^{l_i} \right|}^2 g(z) dz  \\
&\leq& {{\sup_{z \in \jk}} g(z) \over n^2}  \int_{\jk} {\left| \sum_{i=0}^{n-1} z^{l_i} \right|}^2  dz \\
&\leq&  {{\sup_{z \in \jk} } g(z) \over n}.
\end{eqnarray*} 

From this we deduce by the Borel Cantelli  Lemma that $S_{n^2}f / n^2$ converges to zero for almost every $x \in M$. By another use of the Borel Cantelli Lemma we can then interpolate between $n^2$ and ${(n+1)}^2$ showing that for almost every $x \in M$ we have $S_nf(x) / n  \mathop{\longrightarrow} \limits_{n \rightarrow \infty} 0$ which overrules (\ref{eq1}). \carre 

\subsection{} \label{spectral2} Proposition \ref{spectral} has the following immediate corollary
 
\begin{coro*}    Let $(T, M, \mu)$ be a dynamcial system. If for any complex nonzero function $f \in L^2_0 (M, \mu)$,   there  exist $\tau >1$,   a measurable set $E \subset M$ with $\mu (E) > 0$ and a sequence  $k_{n+1} \geq \tau^{n+1} k_n $   such that for every $x \in E$ we have 
\begin{eqnarray} \label{eq2} \limsup_{n\rightarrow \infty} {1 \over [\tau^n]} \left| \sum_{i=0}^{[\tau^n]-1} f(T^{ik_n}x) \right| >0,
\end{eqnarray} 
then the maximal spectral type of the unitary operator associated to $(T,M, \mu)$ is singular.
\end{coro*} 

\noindent {\sl Proof.} The criterion of Proposition \ref{spectral} holds with the set $E$ and the sequence $l_n$ given by: \newline  $k_1, \ldots, [\tau]k_1, \ldots, k_j, 2 k_j, 3 k_j, \ldots, [\tau^j] k_j, k_{j+1}, 2k_{j+1}, \ldots, [\tau^{j+1}] k_{j+1}, \ldots$    \carre

\subsection{} \label{propi} In the sequel we will assume that $(T,M,\mu)$ satisfies $(i)-(iii)$ of Theorem \ref{PASS}. We fix $1 < \tau < \gamma$ and an arbitrary nonzero function $f \in L^2_0(M,\mu)$. For  $\varepsilon > 0$ we define the set 
 $$D_{\varepsilon} = \left\{ x \in M \left|  f(x) \geq 2 \varepsilon  \right. \right\}.$$

Since $f \in L^2_0(M,\mu)$ is not null, there exists $\varepsilon_0>0$ such that
  $\mu(D_{\varepsilon_0}) > 0$. Theorem \ref{PASS} will hold proved if we show that:

\begin{prop*} Under the conditions of Theorem \ref{PASS}, we have that, for $\mu$ a.e. point $x \in D_{\varepsilon_0}$, there exits infinitly many integers $n$ such that 
\begin{eqnarray} \label{eq3} {1 \over [\tau^n]} \sum_{i=0}^{[\tau^n]-1} f(T^{ik_n}x) \geq \varepsilon_0. 
\end{eqnarray} 
\end{prop*}

\subsection{} \label{final}  For $x \in D_{\varepsilon_0}$, let $N(x) \in \RR^+ \cup \lbrace \infty \rbrace$ be such that for every  $n \geq N(x)$, (\ref{eq3}) fails for $x$. The function $N(x)$ is naturally measurable and we have to show that almost surely it is equal to infinity.
This will clearly hold if we prove the following

\begin{prop*} Under the conditions of Theorem \ref{PASS}, 
 for every $N \geq 0$ and for every measurable set $D \subset D_{\varepsilon_0}$, we can find a set $\overline{D} \subset D$ satisfying
\begin{itemize}
 \item $\mu(\overline{D}) > 0$; 

\item For   every $x \in \overline{D}$, (\ref{eq3}) holds for some $n \geq N$. 
\end{itemize}

\end{prop*}


\subsection{} \label{le} Define $f_0 = \min (f, 2 \varepsilon_0)$. To prove Proposition \ref{final} we will need the following Lemma

\begin{lemm*}  There exists $N_0$ such that if $n \geq N_0$ and $B_n$ is a set satisfying $(ii)$ of Theorem \ref{PASS} and $$\int_{B_{n}} f_0(x) d \mu (x) \geq {3 \over 2} \varepsilon_0 \mu(B_{n})$$ 
then there exists  a set $\overline{B}_n \subset B_{n}$ with $\mu(\overline{B}_n) \geq { \mu(B_n) / 5}$    such that (\ref{eq3}) holds for every $x \in \overline{B}_n$.      
\end{lemm*} 
 
\noindent {\sl Proof.} Let $B_n$ and $k_n$ be as in $(ii)$ of Theorem \ref{PASS}. 
For $x \in M$, we use in this proof the notation 
$${S_nf(x)}  := \sum_{i=0}^{[\tau^n]-1} f(T^{ik_n}x).$$

Define  
\begin{eqnarray*} \tilde{B}_n = \bigcup_{i=0}^{[\tau^n]-1} T^{-ik_n} B_n \qquad \hat{B}_n  =  \bigcap_{i=0}^{[\tau^n]-1} T^{-ik_n} B_n. \end{eqnarray*}
Clearly $\hat{B}_n \subset B_n \subset \tilde{B}_n$ and since $\tau < \gamma$, $(ii)$ implies  for $n$  sufficiently large 
\begin{eqnarray} \label{aa} \mu( \tilde{B}_n \triangle \hat{B}_n ) \leq {\varepsilon_0 \over 100} \mu(B_n). \end{eqnarray}

Define $\tilde{f}_0 = f_0$ on $B_n$ and equal to zero otherwise. We then have 

$$\displaystyle{
\int_{\tilde{B}_n}  {S_n \tilde{f}_0(x) \over [\tau^n]} d \mu(x)   = \int_{M}  {S_n\tilde{f}_0 (x) \over [\tau^n]} d \mu(x) = \int_M \tilde{f}_0(x) d \mu(x) = \int_{B_n} {f}_0 d \mu(x)},$$
hence from our hypothesis
\begin{eqnarray} \label{pp} \int_{\tilde{B}_n}  {S_n \tilde{f}_0(x) \over [\tau^n]} d \mu(x) \geq  {3 \over 2} \varepsilon_0 \mu(B_{n}).  
 \end{eqnarray}

On the other hand, since $\tilde{f}_0 \leq 2 \varepsilon_0$ we get

\begin{eqnarray*} \label{qq} \int_{\tilde{B}_n}  {S_n \tilde{f}_0(x) \over [\tau^n]}  d \mu(x)   \leq  \mu(\tilde{B}_n) \varepsilon_0 + \mu \left( \left\{ x \in \tilde{B}_n \left| {S_n \tilde{f}_0 (x) \over [\tau^n]}  \geq \varepsilon_0 \right\} \right. \right) 2 \varepsilon_0
\end{eqnarray*} 
which in light of (\ref{aa}) and (\ref{pp}) leads to 
$$\mu \left(\left\{ x \in \tilde{B}_n \left| {S_n \tilde{f}_0(x) \over [\tau^n]}  \geq \varepsilon_0 \right\} \right. \right) \geq  \left({1 / 4} -{1 / 200} \right) \mu(B_n),$$
which using (\ref{aa}) again  yields
$$\mu \left(\left\{ x \in \hat{B}_n \left|  {S_n \tilde{f}_0(x) \over [\tau^n]}  \geq \varepsilon_0 \right. \right\} \right) \geq {1 / 5}  \mu(B_n),$$ 
which is the desired inequality since $S_n \tilde{f}_0$ and $S_n f_0$ coincide on $\hat{B}_n \subset B_n$.  \carre

\subsection{} { \sl Proof of Proposition \ref{final}.} Let $D$, a measurable subset of $D_{\varepsilon_0}$ such that $\mu(D) >0$, and $N \in \NN$ be fixed. Define $\overline{N}= \sup (N_0,N)$ where $N_0$ is as in Lemma \ref{le}.

 By Vitali's Lemma and properties $(i)$ and $(iii)$, there exists a constant $0< \vartheta <1$ such that, given any ball $B$ in $M$, we can find a  familly of balls $B_{n_i} \subset B$ such that 
\begin{itemize}
\item[(P1)] The $B_{n_i}$ are disjoint;

\item[(P2)] Every $B_{n_i}$ belong to some ${\mathcal C}_n$ with $n \geq \overline{N}$;

\item[(P3)] $\displaystyle{ \mu \left({\bigcup} B_{n_i} \right) \geq \vartheta \mu(B).}$ 
 
\end{itemize}

For $x \in D \subset D_{\varepsilon_0}$, we have $f_0 = 2 \varepsilon_0$. Considering a Lebesgue density point  we obtain, for any $\epsilon >0$,  a ball $B \subset M$ such that 

\begin{itemize}

\item[(B1)] $\displaystyle{ \mu ( B \cap D ) \geq (1-\epsilon) \mu ( B) }$;

\item[(B2)] $\displaystyle{ \int_B f_0 (x) d \mu(x) \geq (2 - \epsilon) \varepsilon_0 \mu(B)}$.

\end{itemize}

We can choose $\epsilon>0 $ arbitrarilly small in (B1), (B2) and  then apply (P1)-(P3) to the above ball $B$. We can hence obtain a ball $B_{n} \in {\mathcal C}_n$ such that $n \geq \overline{N}$ and $\mu( B_n \cap D ) \leq (1- 1 /10) \mu(B_n)$ while $\int_{B_n} f_0 (x) d \mu(x) \geq 3 /2 \varepsilon_0 \mu(B_n)$. We conclude using Lemma \ref{le}. \carre


{

\section{Application: Slow mixing and singular spectrum} \label{S3}

This section is devoted to the proof of Theorem \ref{repara}.

\subsection{Reduction to special flows.} \label{special}

\begin{defi*}(Special flows)
Given a Lebesgue space $L$, a measure preserving transformation $T$ on
$L$ and an integrable strictly positive real function defined on $L$ we define the special
 flow over $T$ and under the {\sl ceiling function} $\varphi$
by inducing on $M({L,T,\varphi})= {{L \times \RR} / \sim }$, where $\sim$ is the
identification $(x, s + \varphi(x)) \sim (T(x),s)$,the action of \begin{eqnarray*}L \times \RR &  \rightarrow &  L \times \RR \\ (x,s) & \rightarrow & (x,s+t). \end{eqnarray*} \end{defi*} 

 If $T$ preserves a unique probability measure $\lambda$, then the special flow will preserve a unique probability measure that is the normalized product measure of $\lambda$ on the base and the Lebesgue measure on the fibers. 

We will be interested in special flows above minimal translations $R_{\a,\a'}$ of the two torus and under smooth functions $
\varphi(x,y) \in C^{\infty}(\TT^2,\RR_+^*)$ that we will denote by $T^t_{\a,\a',\varphi}$. 
For $r \in {\NN} \bigcup \lbrace +\infty \rbrace$, we denote by
$C^r({\TT}^2, \RR)$ the set of real functions on ${\RR}^2$ of class $C^r$
and ${\ZZ}^d$-periodic. We denote by $C^r({\T}^d, \R_{+}^{*})$ the set  of
strictly positive functions in
$C^r({\T}^d, \R)$.
 Without loss of generality, we will consider cealing function $\varphi$ with the property $\int_{\TT^2} \varphi(x,y)dxdy =1$.

In all the sequel we will use the following notation, for $m \in \NN$,
$$S_m \varphi (x,y) = \sum_{l=0}^{m-1} \varphi(x+l\a, y + l \a')$$

With this notation, given $t \in \RR_+$ we have for $z \in \TT^2$
  $$T^t(z,0) = \left( R^{N(t,z)}_{\a,\a'}(z), t -
\varphi_{N(t,z)}(z) \right)$$
where  $N(t,z)$ is the largest integer $m$ such that $t-\varphi_m(x)
  \geq 0$,
 that is the number of fibers covered by  $z$ during its
  motion under the action of the flow until time $t$.

By the equivalence between special flows and reparametrizations Theorem \ref{repara} follows if we prove
\begin{theo*} There exists a vector $(\a,\a') \in \RR^2$ and $\varphi \in C^\infty(\TT^2, \RR_+^*)$ 
such that  the special flow $T^t_{\a,\a',\varphi}$  is mixing and satisfies $(i)-(iii)$ of Theorem \ref{PASS}, which implies that the 
spectral type of the flow is purely singular. 
\end{theo*} 

The equivalence between the above theorem and Theorem \ref{repara} is  standard  and can be found in \cite{mixing}, Section 4.

We will now undertake the construction of the special flow $T^t_{\a,\a',\varphi}$. We will first choose a special translation vector on $\TT^2$, then we will give two criteria on the Birkhoff sums of the special function $\varphi$ above $R_{\a,\a'}$ that will guarantee mixing and SCPA. Finally, we build a smooth function $\varphi$ satisfying these criteria. 

\subsection{Choice of the translation on $\TT^2$.} \label{R}

Given a real number $u$, we will use the following notations: $[u]$ to indicate the integer part of $u$, $\lbrace u \rbrace$ its fractional part and $\tn u \tn$ its closest distance to integers.   Let $\a$ be
an irrational real number, then  there exists a sequence of  rationals 
${\lbrace {p_n \over
  q_n} \rbrace}_{n \in \NN} $, called the convergents of $\a$, such that
\begin{eqnarray}  \label{best} \tn q_{n-1} \a \tn < \tn k \a \tn, \ \ \forall k < q_{n} \end{eqnarray} 
and for any $n \in \NN$ 

\begin{eqnarray} \label{equation}
  { 1  \over q_n (q_n + q_{n+1})} \leq {(-1)}^n ( \a - {p_n \over q_n}) \leq   { 1  \over q_n q_{n+1}}.
\end{eqnarray}

We recall also that any
irrational number $\a \in \RR - \QQ$ has a writing in continued fraction
$$ \a = [a_0,a_1,a_2,...]= a_0 + 1 / (a_1 + 1 / (a_2+...)),$$               
where ${\lbrace a_i \rbrace}_{i \geq 1}$ is a sequence of integers
$\geq 1$,
$a_0 = [\a]$. Conversely,  any sequence  ${\lbrace a_i \rbrace}_{i \in \NN}$
 corresponds to a unique number $\a$. The convergents of $\a$ are
given by the $a_i$ in the following way:
\begin{eqnarray*}
p_n = a_n p_{n-1} + p_{n-2} \ \ &{\rm for}& \ \ n \geq 2, \ \  p_0=a_0, \ \ p_1=a_0a_1 + 1, \\
q_n = a_n q_{n-1} + q_{n-2} \ \ &{\rm for}& \ \ n \geq 2, \ \ q_0=1,
 \ \ q_1=a_1. 
\end{eqnarray*}

Following \cite{Y} and as in \cite{mixing}, we take $\a$ and $\a'$ satisfying 
\begin{eqnarray}
q'_{n} &\geq& e^{3q_n}, \label{111} \\
q_{n+1} &\geq& e^{3q'_n}. \label{222}
\end{eqnarray}

Vectors $(\a,\a') \in \RR^2$ satisfying (\ref{111}) and (\ref{222}) are obtained by an adequate inductive choice of the sequences  $a_n(\a)$ and $a_n(\a')$. Moreover, it is easy to see that the set of vectors $(\a,\a') \in \RR^2$ satisfying (\ref{111}) and (\ref{222}) is a continuum (Cf. \cite{Y}, Appendix 1).

\subsection{Mixing criterion} \label{criterion mixing}
 
We will use the criterion on mixing for a special flow $T^t_{\a,\a',\varphi}$ studied in \cite{mixing}. It is based on the uniform stretch of the Birkhoff sums $S_m \varphi$ of the ceiling function above the $x$ or the $y$ direction alternatively depending on whether $m$ is far from  $q_n$ or from $q'_n$. From \cite{mixing}, Propositions 3.3, 3.4 and 3.5 we have the following sufficient
 mixing criterion:
 
\begin{prop*}[Mixing  Criterion ] Let $(\a,\a')$ be as in (\ref{111}) and (\ref{222}) and $\varphi \in C^2(\TT^2, \RR^*_+)$.
If for every $n \in \NN$ sufficiently large, we have a two sets $I_n $ and $I_n'$, each one  being equal to the circle minus two intervals whose lengths converge to zero, such that:

\begin{itemize} 

\item $\displaystyle{ m \in \left[{e^{2q_n} / 2} , 2 e^{2q'_n} \right] \ \ \  \Longrightarrow  \left| D_x S_m \varphi (x,y) \right| \geq {m \over e^{q_n}} }$, for any $y \in \TT$ and any $x$ such that $\lbrace q_n x \rbrace \in I_n$;
 
\vspace{0.1cm} 

\item $\displaystyle{ m \in \left[{e^{2q'_n} / 2} , 2 e^{2q_{n+1}} \right]  \ \ \Longrightarrow  \left| D_y S_m \varphi (x,y) \right| \geq {m \over e^{q'_n}} }$, for any $x \in \TT$ and any $y$ such that $\lbrace q'_n y \rbrace \in I_n'$;

\end{itemize}

Then the special flow $T^t_{\a,\a',\varphi}$ is mixing.

\end{prop*}

\subsection{Criterion for the existence of slowly coalescent periodic approximations} \label{criterion approximations}

 We give now a condition on the Birkhoff sums of $\varphi$ above $R_{\a,\a'}$ that is sufficient to insure SCPA for $T^t_{\a,\a',\varphi}$ on $M= M(\TT^2, R_{\a,\a'}, \varphi)$:
\begin{prop*} If for $n$ sufficiently large, we have for any $x$ such that $1 / n^2 \leq \lbrace q_n x \rbrace \leq 1 /n - 1 / n^2$  and for any $y \in \TT$ 
\begin{eqnarray} \label{r1} \left| S_{q_nq'_n} \varphi(x,y) - q_nq'_n  \right| \leq {1 \over e^{q_n}},
\end{eqnarray}
then the special flow $T^t_{\a,\a',\varphi}$ has  { slowly coalescent periodic approximations} as in Definition \ref{PASS}.

\end{prop*}

\noindent {\sl Proof.} Let $C_n$ be the set of points  $(x,y,s) \in M$ satisfying ${2 / n^2} \leq \lbrace q_n x \rbrace \leq {1 / n} - {2 / n^2}$. It follows from the definiton of special flows and (\ref{r1}) that for $(x,y,s) \in M$ such that  $1 / n^2 \leq \lbrace q_n x \rbrace \leq {1/ n} - 1 / n^2$ we have 
$$T^{q_nq'_n}(x,y,s) = \left(x+q_nq'_n \a, y+ q_nq'_n \a', s +  S_{q_nq'_n} \varphi(x,y) - q_nq'_n \right)$$
but from (\ref{equation}) we have that  $\tn q_nq'_n \a \tn \leq q'_n / q_{n+1} = o(e^{-q_n})$ as well as  
  $\tn q_nq'_n \a' \tn \leq q_n / q'_{n+1} = o(e^{-q_n})$. Therefore (\ref{r1})  implies that $d(T^{q_nq'_n}(x,y,s), (x,y,s)) \leq {2 / e^{q_n}}$. It is therefore possible to cover $C_n$ with a collection of balls ${\mathcal C}_n$ such that each ball $B \in {\mathcal C}_n$ has radius less than $1 / nq_n$ and satisfies 
 $ \mu \left( T^{q_nq'_n} B \triangle B \right) \leq e^{-n} \mu(B)$ which yields  conditions $(i)$ and $(ii)$ of Definition \ref{PASS}.


On the other hand it is clear from the difference of scale between the successive terms of the sequence  $q_n$ that the sets $C_n$ are almost independent and the fact that $\mu(C_n) \geq {1 / n} \inf_{(x,y) \in \TT^2} \varphi(x,y)$ then implies by the Borel Cantelli Lemma that $\mu  \left( \mathop{\bigcap} \limits_{m \in \NN} \mathop{\bigcup} \limits_{n \geq m} {C}_n \right) = 1$, which is condition (iii) of the Definition \ref{PASS}. \carre

\subsection{Choice of the ceiling function $\varphi$.}  \label{xtilde}

Let $(\a,\a')$ be as above and define 
$$ f(x,y) = 1+ \sum_{n \geq 2} X_n(x) + Y_n(y)$$ 
where 
\begin{eqnarray}
\label{XX} X_n(x) &=& {1 \over e^{q_n}} \cos (2\pi q_nx) \\
\label{YY} Y_n(y) &=& {1 \over e^{q'_n}} \cos (2 \pi q'_n y).
\end{eqnarray}

Relying on the Proposition-Criterion \ref{criterion mixing} stated above, we proved in \cite{mixing} that the flow $T^t_{\a,\a',f}$ is mixing. In order to keep this criterion valid but have in addition the conditions of Criterion \ref{criterion approximations} satisfied we modify the ceiling function in the following way:

$\bullet $ We keep $Y_n(y)$ unchanged.

\vspace{0.1cm} 

$\bullet$ We replace $X_n(x)$ by a trigonometric polynomial $\tilde{X}_n$ with integral zero, that is essentially equal to $0$ for $\lbrace q_nx \rbrace < 1 /n$ and  whose  derivative has its absolute value  bounded from below by ${1 / e^{q_n}}$ for 
 $ \lbrace q_n x \rbrace \in [2/n, 1/2 - 1/n] \cup [1/2+2/n, 1-1/n]$.
The first listed  properties of $\tilde{X}_n$ will yield Criterion \ref{criterion approximations} while the  lower bound on the absolute value of its derivative will insure Criterion \ref{criterion mixing}.

More precisely, the following Proposition enumerates some properties that we will require on $\tilde{X}_n$ and its Birkhoff sums that will be sufficient for our purposes, and that we will realize with a specific construction at the end of the section.

\begin{prop*}  Let $(\a,\a')$ be as in Section \ref{R}. There exists a sequence of trigonometric polynomials $\tilde{X}_n(x)$ satisfying

\begin{enumerate}

\item[(1)] \label{x1}  $\displaystyle{ \int_{\TT} \tilde{X}_n(x) dx=0}$;

\item[(2)]  \label{x2} For any $r \in \NN$, there exists $N(r) \in \NN$ such that for every $n \geq N(r)$, $\displaystyle{ {\parallel \tilde{X}_n \parallel}_{C^r} \leq {1 \over e^{q_n\over 2}}}$;

\item[(3)]  \label{x3} For $\displaystyle{ \lbrace q_n x \rbrace \in [0, {1 / n} ]}$, 
 $\displaystyle{ |\tilde{X}_n(x)| \leq {1 \over {q'_n}^2 } }$;
  
\item[(4)]  \label{x4} For $\displaystyle{ \lbrace q_n x \rbrace  \in [2 /n,  {1 / 2} - {1 / n}]}$,   $\displaystyle{\tilde{X}_n'(x)  \geq {2 \over e^{q_n}}}$,  as well as 
\newline  for  $\displaystyle{ \lbrace q_n x \rbrace  \in [1/2 + 2/n, 1 - {1 / n}]}$,   $\displaystyle{\tilde{X}_n'(x)  \leq - {2 \over e^{q_n}}}$;

\item[(5)]  \label{x5} For $n \in \NN$ sufficiently large,  $\displaystyle{\left\| S_{q_n} \sum_{l\leq n-1} \tilde{X}_{l} \right\| \leq {1 \over {q'_n}^2}}$; 
 
\item[(6)]   \label{x6} For $n \in \NN$ sufficiently large, we have for any $m \in \NN$, $\displaystyle{\left\| S_{m} \sum_{l\leq n-1} \tilde{X}'_{l} \right\| \leq q_n}$.


\end{enumerate}

\end{prop*}

Before we prove this proposition, let   us show how it allows to produce the example of Theorem \ref{special}.

\subsection{Proof of Theorem \ref{special}.}  Define for some $n_0 \in \NN$

\begin{eqnarray} \label{varphi} 
\varphi(x,y) = 1+ \sum_{n=n_0}^{\infty} \tilde{X}_n(x) + Y_n(y)
\end{eqnarray}
where $Y_n$ is as in (\ref{YY}) and $\tilde{X}_n$ is as in the proposition above. From (\ref{YY}) and Property (2) of $\tilde{X}_n$, we have that $\varphi \in C^\infty (\TT,\RR)$.
Also from  (\ref{YY}) and (2) again, we can choose $n_0$ sufficiently large so that $\varphi$ is strictly positive.  We then have

\begin{theo*} Let $(\a,\a') \in \RR^2$ be as in Section \ref{R} and $\varphi$ be given by (\ref{varphi}). Then the special flow $T^t_{\a,\a',\varphi}$ satisfies the conditions of Propositions \ref{criterion mixing} and \ref{criterion approximations} and is therefore mixing with a singular maximal spectral type. 
\end{theo*}

\noindent {\sl Proof.}  The second part of Proposition \ref{criterion mixing} is valid exactly as in \cite{mixing} since $Y_n$ has not been modified. Briefly, the reason is that due to 
(\ref{equation}) and (\ref{111})-(\ref{222}) we have  
$Y_n(y+l\a') \sim Y_n (y)$ for every $l \leq m  \ll q'_{n+1}$ so that $|S_m Y_n'|$ is large as required for $m \in [e^{2q'_n}/2, 2 e^{q_{n+1}}]$. Meanwhile,    $S_m \sum_{k < n} Y_k$ is much smaller because these lower frequencies behave as controlled coboundaries for this range of $m$. As for $S_m \sum_{k>n} Y_k $, it is still very small since $m \ll e^{q'_{n+1}}$.
The latter phenomena will be further explicited and used in the sequel. 

Let $m \in [e^{2q_n} / 2, 2e^{2q'_n}]$ and define $I_n := [3/n, 1/2  - 2 / n] \cup [1/2 + 3/n,1-2/n]$. For $x$ such that $\lbrace q_nx \rbrace \in I_n$,  it follows from (\ref{equation}) that for any $l \leq m$, $2 / n \leq  \lbrace q_n (x+m\a) \rbrace \leq 1 / 2  - 1 / n$. Hence, by Property (4) of $\tilde{X}_n$ 
$$S_m \tilde{X}_n'(x) \geq {2m  \over e^{q_n}}.$$
 On the other hand, Properties (2) and (6) imply that 
\begin{eqnarray*} \left\| S_m \varphi'- S_m \tilde{X}_n'  \right\| &\leq& 
  \left\| S_m \sum_{l<n} \tilde{X}_l' \right\| +
  \left\| S_m \sum_{l>n} \tilde{X}_l' \right\| \\
&\leq& q_n + m \sum_{l \geq n+1} {1 \over e^{{q_l \over 2}}} \\
&\leq& q_n + {2m \over e^{{q_{n+1} \over 2}}} \\
&=& o({m \over e^{q_n}}) \end{eqnarray*}
for the current range of $m$. With an exactly  similar computation for the other part of $I_n$, the criterion of Proposition \ref{criterion mixing} holds  proved.

Let now $x $ be as in Proposition \ref{criterion approximations}, that is $1 / n^2 \leq \lbrace q_n x \rbrace \leq 1 /n - 1 / n^2$. From (\ref{R}) we have for any $l \leq q_n q'_n$ that $0 \leq \lbrace q_n(x+l\a) \rbrace \leq 1 / n$, hence  Property (3) implies 
\begin{eqnarray} |S_{q_nq'_n} \tilde{X}_n (x) | \leq {q_n \over q'_n} \label{101} \end{eqnarray}
the latter being very small compared to $ 1 / e^{q_n}$ since $q'_n \geq e^{3q_n}$. From Properties (5) and (2)  we get for $n$ sufficently large
\begin{eqnarray} \parallel S_{q_nq'_n} \sum_{l \neq n}  S_m \tilde{X}_l   \parallel &\leq& {1 \over q'_n} + q_nq'_n \sum_{l \geq n+1} {1 \over e^{{q_l \over 2}}} \nonumber  \\
&\leq&   {2 \over q'_n}.  \label{102}
\end{eqnarray}

On the other hand, it follows from (\ref{best}) and (\ref{equation}) that for any $y \in \TT$, for any $|j| < q'_{n}$, we have 
\begin{eqnarray}
|S_{q'_n} e^{i2 \pi jy}| &=& \left| \sin(\pi j q_n' \a') \over \sin(\pi j \a') \right|  \nonumber \\
&\leq&  {\pi j q_n' \over q'_{n+1}}, \label{103} 
\end{eqnarray}
which yields for $Y_l$ as in (\ref{YY})
\begin{eqnarray} \left\| S_{q'_n} \sum_{l<n} Y_l \right\| = o({1 \over e^{q'_n}}) \label{104} \end{eqnarray} 
while clearly 
\begin{eqnarray} \left\| S_{q'_n} \sum_{l \geq n} Y_l \right\| = o({ 1 \over e^{{q'_n \over 2} } }). \label{105} \end{eqnarray}

In conclusion, (\ref{r1}) follows from (\ref{101}), (\ref{102}), (\ref{104}) and (\ref{105}). \carre

\medskip 

It  remains to construct $\tilde{X}_n$ satisfying (1)-(6).

\subsection{Proof of Proposition \ref{xtilde}.}

Consider on $\RR$ a $C^\infty$ function, $0 \leq \theta \leq 1$ such that 
\begin{eqnarray*}
\theta(x) &=& 0 \ {\rm for } \ x \leq 1,  \\
\theta(x) &=&  1  \ {\rm for } \ x \geq 2.
\end{eqnarray*}

Then, for $n \in \NN$, define on $\RR$ the $C^\infty$ function

$$\xi_n(x) = \int_0^x \left[ \theta(nq_nt) - \theta \left( nq_n (t - {1 \over 2q_n} + {2 \over nq_n}) \right) \right] dt.$$

It is easy to check the following 
 
\begin{itemize}
\item $\displaystyle{ \xi_n(x) =0 \qquad {\rm for} \ x \leq {1 \over nq_n} }$; 
\item $\displaystyle{ \xi_n'(x) = 1 \qquad {\rm for} \ x \in \left[ {2 \over nq_n}, {1 \over 2q_n} - {1 \over nq_n} \right]}$;

\item $\displaystyle{ 
 \xi_n(x) = \xi_n({1 / 2}) \qquad {\rm for} \ x \geq {1 \over 2q_n} }$.

\end{itemize}

We then introduce the function 
$$\varsigma_n(x) = \xi_n(x) - \xi_n({1\over 2}) \theta \left[ nq_n \left(x -{1 \over 2q_n} + {2 \over nq_n} \right) \right]$$
and define for $x \in [0,1/q_n]$ the function
$$\hat{X}_n(x) = {3 \over e^{q_n}} \left( \varsigma_n(x) - \varsigma_n (x-{1 /2}) \right)$$
that we extend to a $C^\infty$ function over the circle periodic with period $1 / q_n$. It satisfies

\begin{itemize}

\item $\displaystyle{ \int_{\TT} \hat{X}_n (x) dx = 0}$;

\item For any $r \in \NN$, \ \ \ for $n$ sufficiently large   $\displaystyle{ {\left\| \hat{X}_n \right\| }_{C^r} \leq {1 \over e^{3q_n \over 4}}}$;
   
\item $\displaystyle{ \hat{X}_n(x) =0 \qquad {\rm for} \ \lbrace q_n x \rbrace  \in [0, {1 / n}] }$; 
 
\item $\displaystyle{ \hat{X}_n'(x) = {3 \over e^{q_n}}  \qquad {\rm for} \ \lbrace q_n x \rbrace  \in \left[ {2 / n}, {1 / 2} - {1 / n} \right]}$, \ and 
\newline $\displaystyle{ \hat{X}_n'(x) = {-3  \over e^{q_n}}  \qquad {\rm for} \ \lbrace q_n x \rbrace  \in \left[ {1 / 2} + {2 / n}, 1  - {1 / n} \right]}.$
  \end{itemize}

Finally, we  consider the Fourier series of $\hat{X}_n(x) = \sum_{k \in \ZZ}  \hat{X}_{n,k} e^{i2\pi k x}$
and let

$$\tilde{X}_n := \sum^{q_{n+1}-1}_{k=-q_{n+1}+1} \hat{X}_{n,k} e^{i 2 \pi k x}.$$ 

The Fourier coefficients $f_k$ of  a function $f \in C^\infty(\TT,\RR)$ satisfy for any $k \in \ZZ$
\begin{eqnarray} {(2 \pi)}^{r-1} {|k|}^r |f_k| \leq {\| f \|}_{C^r} \leq \sup_{k \in \NN} {(2\pi |k|)}^{r+2} |f_k|. \label{four} \end{eqnarray}

Hence, we have for any $r \in \NN$
\begin{eqnarray*} 
{\left\| \tilde{X}_n - \hat{X}_n \right\|}_{C^r} &\leq&  \sum_{|k| \geq q_{n+1}} {(2\pi k)}^{r} |\hat{X}_{n,k}| \\
&\leq& {1 \over 2 \pi} {\| \hat{X}_n \|}_{C^{r+2}} \sum_{|k| \geq q_{n+1}} {1 \over k^2} \\
&=& o({1 \over {q'_n}^2})
\end{eqnarray*}
which allows to check (1), (2), (3) and (4) for $\tilde{X}_n$ from the properties of $\hat{X}_n$.

\medskip

\noindent  {\sl Proof of Properties (5) and (6).} We have 
due to our truncation
\begin{eqnarray} \label{sub} \tilde{X}_n(x) = \psi_n(x+\a) - \psi_n(x) \end{eqnarray}
where 
$$\psi_n(x) = 
\sum^{q_{n+1}-1}_{k=-q_{n+1}+1} \psi_{n,k}  e^{i 2 \pi k x}$$
with
$$\psi_{n,0}=0 \qquad {\rm and \ for \ k \neq 0, \ \ } \psi_{n,k} = {\hat{X}_{n,k} \over  e^{i2\pi k \a} -1}.$$

Since $|k| < q_{n+1}$, it follows from (\ref{best})  that 
$$|\psi_{n,k} | \leq q_{n+1} |\hat{X}_{n,k}|$$
which with (\ref{four}) implies 
\begin{eqnarray*} {\| \psi_n \|}_{C^r} &\leq& 2 \pi q_{n+1} {\| \hat{X}_n \|}_{C^{r+2}} \\
&\leq&  2 \pi {q_{n+1} \over e^{{3 q_n \over 4}}} 
\end{eqnarray*}
for sufficiently large $n$. Hence, from (\ref{sub}) and (\ref{equation}) we get
\begin{eqnarray*}
\left\| S_{q_n} \sum_{l \leq n-1} \tilde{X}_l \right\| &\leq& {1 \over q_{n+1}}   \sum_{l \leq n-1} {\| \psi_l \|}_{C^1} \\ 
&\leq& {1 \over q_{n+1}} \sum_{l \leq n-1}   {q_{l+1} \over e^{{3 q_l \over 4}}} \\
&\leq& {q_n \over q_{n+1}} 
\end{eqnarray*} 
so that property (5) follows. Similarly, we have for sufficiently large $n$ 
\begin{eqnarray*}
\left\| S_{m} \sum_{l \leq n-1} \tilde{X}_l' \right\| &\leq& 2 \sum_{l \leq n-1} {\| \psi_l \|}_{C^1} \\ 
&\leq& q_n. 
\end{eqnarray*}
\carre


\frenchspacing
\bibliographystyle{plain}
 
\end{document}